\documentclass[12pt]{amsart}
\usepackage{amssymb,amsmath,amsthm}

\begin{document}

\theoremstyle{plain}
\newtheorem{theorem}{Theorem}
\newtheorem{lemma}[theorem]{Lemma}
\newtheorem{proposition}[theorem]{Proposition}
\newtheorem{corollary}[theorem]{Corollary}
\newtheorem{conjecture}[theorem]{Conjecture}
\newtheorem*{main}{Main Theorem}

\def\mod#1{{\ifmmode\text{\rm\ (mod~$#1$)}
\else\discretionary{}{}{\hbox{ }}\rm(mod~$#1$)\fi}}

\theoremstyle{definition}
\newtheorem*{definition}{Definition}

\theoremstyle{remark}
\newtheorem*{example}{Example}
\newtheorem*{remark}{Remark}
\newtheorem*{remarks}{Remarks}

\newcommand{\qq}{{\mathbb Q}}
\newcommand{\rr}{{\mathbb R}}
\newcommand{\nn}{{\mathbb N}}
\newcommand{\zz}{{\mathbb Z}}

\newcommand{\al}{\alpha}
\newcommand{\be}{\beta}
\newcommand{\ga}{\gamma}

\newcommand{\ep}{\epsilon}
\newcommand{\la}{\lambda}
\newcommand{\de}{\delta}
\newcommand{\Del}{\Delta}

\thanks{This material is based in part upon work of the author,
  supported by the USAF under DARPA/AFOSR MURI Award
  F49620-02-1-0325. Any opinions, findings, and conclusions or
  recommendations expressed in this publication are those of the
  author and do not necessarily reflect the views of these agencies.}

\title{Clean lattice tetrahedra}

\author{Bruce Reznick}
\address{Department of Mathematics, University of 
Illinois at Urbana-Champaign, Urbana, IL 61801} 
\email{reznick@math.uiuc.edu}

\begin{abstract}
A clean lattice tetrahedron is a non-degenerate tetrahedron with the property
that the only lattice points on its boundary are its vertices. We
present some new proofs of old results and some new results on clean
lattice tetrahedra, with an emphasis on counting the number of its interior
lattice points and on computing its lattice width.
\end{abstract}
\date{\today}

\maketitle

\section{Introduction and Overview}

Let $T = T(v_1,\dots,v_n) = conv(v_1,\dots,v_n)$ be a non-degenerate
simplex with 
vertices $v_j \in \zz^n$. We say that $T$ is \emph{clean} if
there are no non-vertex lattice points on the boundary of $T$. Let
$i(T) = \# \{ int(T) \cap \zz^n \}$ denote the number of lattice points
in the interior of a clean lattice simplex $T$. If $i(T) = k$, then 
$T$ is called a \emph{$k$-point lattice simplex}. If $i(T) = 0$,
then $T$ is called  \emph{empty}. This paper is mainly concerned with
clean tetrahedra.

Pick's Theorem says that the area of a clean lattice triangle $T$ is equal to
$i(T) + 1/2$. Reeve \cite{ree}
showed in 1957 that
there are empty lattice tetrahedra having arbitrarily large volume. By
constrast, if $T$ is a  (not necessarily clean) lattice tetrahedron
with $k \ge 1$ interior points, then Hensley \cite{he} showed in 1983 that
there is an upper bound on the volume of $T$ depending on $k$. 
Any lattice tetrahedron determines an (affine) lattice $\Lambda$,
and if $|\zz^3/\Lambda| = m$; that is, if there are $m$ lattice points
in the fundamental parallelepiped, then the volume of 
 that parallelepiped equals $m$.  The volume of the corresponding
tetrahedron is then equal to $m/6$, but there seems to be no easy way to
determine the number of lattice points it contains.

In this paper, we give a unified discussion of clean lattice
tetrahedra. We begin with preliminaries in section two.
Two tetrahedra $T$ and $T'$ are {\it equivalent} 
if there is an affine unimodular map which takes the vertices of $T$
into the vertices of $T'$ in some order.
 For $(a,b,n) \in \zz^3$, we define the
tetrahedron $T_{a,b,n}$, which has vertices $(0,0,0)$, $(1,0,0)$,
$(0,1,0)$ and $(a,b,n)$, and we give necessary and sufficient
conditions under which $T_{a,b,n}$ and $T_{a',b',n'}$ are
equivalent. A crucial ``hidden'' parameter is $c = 1 - a - b$.

In section three, we show that every clean lattice tetrahedron is
equivalent to some $T_{a,b,n}$, where $\gcd(a,n) = \gcd(b,n) =
\gcd(c,n) = 1$ and $0 \le a,b \le n-1$. (Reeve had originally
discussed these conditions in the context of
empty lattice tetrahedra.) We then review our 1986 result
\cite{rz} that
$$
i(T_{a,b,n}) = \#\left\{t: 1\le t \le n-1\quad \text{and} \quad
  \{\tfrac{t(a+b-1)}n \} 
+ \{ \tfrac{-ta}n \} + \{ \tfrac{-tb}n \} +
\{ \tfrac tn \} = 1\right\}.
$$
Using this, we give a
new and shorter proof of White's 1964 theorem \cite{w} that a lattice
tetrahedron is empty if and only if it is equivalent to $T_{0,0,0}$ or
some $T_{1,b,n}$, where $n \ge 2$, $1 \le b \le n-1$ and $\gcd(b,n) =
1$. We also discuss bounds on $i(T_{a,b,n})$. It is not hard to show
that $i(T_{3k,3k,3k+1}) = k$. Han Duong
has proved that $i(T_{2k+1,4k+3,12k+8}) = k$, and used the formula
above to conjecture that for all clean tetrahedra,
$$
\frac{n-1}3 \ge i(T_{a,b,n}) \ge \frac{n-8}{12},
$$
with the extreme examples given by the two aforementioned families.

Section four is devoted to 1-point lattice tetrahedra. Suppose $T$ is
such a tetrahedron, with interior point $w$. We give a new proof of
our earlier result that there are only seven possible sets of
barycentric coordinates for $w$ with respect to the vertices of $T$.
If two 1-point lattice tetrahedra are equivalent, then their interior
points have the same barycentric coordinates, but the converse is
false: $T_{3,3,4}$ and $T_{3,7,20}$ have different volumes and so are
not equivalent, but each has a single interior lattice point at the
centroid. Mazur has \cite{m} recently showed that, up to
equivalence, these are the only two such 1-point lattice
tetrahedra. (This was the fruit of an undergraduate research project.)
We show that in the other six cases of barycentric
coordinates, there is exactly one equivalence class of 1-point
tetrahedra. After an early version of this paper was distributed,
Julian Pfeifle pointed out that this result had been proved
recently by A. Kasprzyk \cite{ka} . The proof here seems sufficiently
different to merit publication. Kasprzyk's paper is motivated by a
question in toric varieties on the classification of toric Fano
3-folds with terminal singularities. In his discussion, the 1-point
tetrahedra are arranged so that the interior point is the origin. 

Finally, in section five, we discuss the lattice width of clean
tetrahedra. White's Theorem showed that an empty tetrahedron lies in two
consecutive planes of lattice points; that is, an empty tetrahedron
has lattice width one. We  show that each 1-point tetrahedron
lies in three consecutive
planes of lattice points, and so has lattice width two. 
Since $T_{3k,3k,3k+1}$ also has lattice width two, there is no
deterministic connection between $i(T)$ and its lattice width; however,
we conjecture that the lattice width of a clean tetrahedron $T$ is
bounded above by $i(T) + 1$. We also show that the lattice width of
$T$ is $\mathcal O(n^{1/3})$ and that $T_{n,n^2,n^3-1}$ has lattice
width $n$, so this bound is asymptotically best possible, up to
multiplicative constant.

Most of the literature is not fastidious about the
existence of lattice points on the boundary of a simplex, unless there are no
interior points.  A principal result of
\cite{rz} was that, if $S \in \rr^n$ is a clean lattice simplex with
exactly $k$ interior points, then there is an upper bound, depending on
$k$ and $n$, for the
denominators of the barycentric coordinates of these points. This
result was subsumed by the stronger and 
essentially simultaneous work of D.~Hensley
\cite{he}, who proved that if $S \in \rr^n$ is a lattice simplex
with $k \ge 1$ interior points, then there are bounds on the volume of $S$.
These bounds were subsequently improved by J.~Lagarias
and G.~Ziegler \cite{lz} and by O.~Pikhurko \cite{p}. A special case of 
Pikhurko's bound
shows that the volume of a lattice tetrahedron with one interior
point is $\le \frac{31^3}{3!352}< \frac {85}6$. This result could be
combined with Theorem 4(i) below and some computer searching to
determine all 1-point lattice tetrahedra up to equivalence. Nevertheless, we
believe it is worthwhile to give a proof in which all computations are
explicitly presented.
Note that the tetrahedron with vertices at the origin and
$\{4e_j\}$, $1 \le j \le 3$, has volume $\frac{64}6$  and a single
interior point $(1,1,1)$, 
along with many boundary points. It is plausible to believe that this
volume is maximal among such tetrahedra. 

The author would like to
thank his fellow organizers of the 2003 Snowbird Conference on Integer
points in Polyhedra -- Sasha Barvinok, Matthias Beck, Christian Haase,
Mich\'ele Vergne, and Volkmar Welker -- for the invitation to join
them in that enterprise, an experience which revived my interest in this
subject. 
Jeff Lagarias reminded me there of the intuitively contradictory
results that empty simplices have unbounded volume, but 1-point
simplices do not. That conversation motivated a
short contribution to the problems article  \cite{bcfhkrrs} from the
Snowbird conference, which has grown into the present paper.

The author would
also like to thank Julian Pfeifle (for pointing out \cite{ka}) and
Alex Kasprzyk (for his insights on toric geometry) and  his
students Han Duong, Ricardo Rojas and Melissa Simmons (for
their patience in listening to earlier versions of this work during
various seminars in the summer of 2004.)

\section{Preliminaries}
Let $T = T(v_1,v_2,v_3,v_4) = conv(v_1,v_2,v_3,v_4)$ be a non-degenerate
tetrahedron in 
$\rr^3$. Every point $w \in \rr^3$ has a unique set of \emph
{barycentric coordinates} with respect to $T$; namely $\la_j:=
\la_{j,T}(w) \in \rr$, $1 \le j \le 4$, so that 
\begin{equation}\label{E:BC}
w = \sum_{j=1}^4 \la_jv_j,\qquad \sum_{j=1}^4 \la_j = 1.
\end{equation}
If \eqref{E:BC} holds, we write $BC_T(w):=
(\la_1(w),\la_2(w),\la_3(w),\la_4(w))$. 
 If the vertices of $T$ are permuted, $T$ as a geometric
object is unchanged, but the coordinates of $BC_T(w)$ are permuted.

Observe that $w \in T$ if and only if $\la_j(w) \ge 0$ for all $j$ and 
$w \in int(T)$
 if and only if  $\la_j(w) > 0$ for all $j$.
If $T$ is clean, $w \in T$ and $\la_j(w) = 0$ for some $j$, 
 then $w=v_k$ for some $k$ and $BC_T(w)$ is a unit vector.

Recall that $x \in \rr$ can be written $x = \lfloor x \rfloor +
\{x \}$, where $\lfloor x \rfloor \in \zz$ and $\{x \} \in [0,1)$. 
If $\sum x_j \in \zz$, then so is $\sum \{x_j\}$. In particular, 
 if $x,y \in \rr$ and $x+y \in \zz$, but  $x,y
\notin \zz$, then $\{x\} + \{y\} = 1$; thus, if $x \notin \zz, \lfloor
-x \rfloor + \lfloor x \rfloor = -1$. Further, if $m$ is an integer
and $x+m \in [0,1)$, then $m = - \lfloor x \rfloor$ and $x+m = \{x\}$.
If $a$ and $n$ are integers, then $a \equiv n\{\frac an \} \mod n$.

Let $\mathcal L$ denote the set of affine unimodular maps 
$f:\zz^3 \to \zz^3$
given by $f(v) = Mv + u$, where $M \in M_3(\zz)$, $\det(M) = \pm 1$
and $v \in \zz^3$. Then $f^{-1} \in \mathcal L$ as well and $f$ is  an
bijection of $\zz^3$ to itself. Since
$$
\sum_{j=1}^4 \la_jf(v_j) = \sum_{j=1}^4 \la_j(Mv_j +u) = M\left(
\sum_{j=1}^4 \la_jv_j\right) + \left( \sum_{j=1}^4 \la_j\right) u =
f\left(\sum_{j=1}^4 \la_jv_j\right),
$$ 
$f$ preserves barycentric coordinates; thus,  $f(T) \cap \zz^3 = f(T
\cap \zz^3)$, with boundary and interior points mapped to boundary and
interior points.
For this reason, it makes sense to classify lattice tetrahedra up to
the action of $\mathcal L$. Following \cite{rz}, given lattice
tetrahedra $T = T(v_j)$ and $T' = T(v_j')$, we say that $T$ and $T'$
are \emph {equivalent} ($T \approx T'$) if there exists $f\in \mathcal
L$ so that 
$\{v_j'\} = \{f(v_j)\}$.  It is not necessary that $f$ preserve the
order of the vertices.

The class $\mathcal L$ contains translations and reflections, of
course. It also contains shears; of particular interest is the map
$(x,y,z) \mapsto (x-mz,y-nz,z)$ for $m,n \in \zz$, where $m$  and $
n$ are chosen by the Euclidean algorithm  so that $0 \le x-mz, y- nz <
|z|$. We call this  a \emph{Euclidean shear}.  
If $r,s \in \zz$ and $g =
\gcd(r,s)$, then there exist  $r',s',m,n \in \zz$ so that $r = gr', s
= gs'$, and  $mr' + ns' = 1$. The map sending
 $(x_j,x_k)$ to $(mx_j +
nx_k,-s'x_j+r'x_k)$, and fixing the other coordinate, has determinant $mr' +
ns' = 1$ and sends $(r,s)$ to $(g,0)$. We  call this a \emph {tweak}. 

For $(a,b,n) \in \zz^3$, $n \neq 0$, we define a standard family of
tetrahedra: 
$$
T_{a,b,n} := T( (0,0,0), (1,0,0), (0,1,0), (a,b,n) ).
$$
The face containing $(0,0,0), (1,0,0), (0,1,0)$ is the \emph{base} of
$T_{a,b,n}$ and  $(a,b,n)$ is the \emph{top}.
 The reflection $(x,y,z) \mapsto (x,y,-z)$ shows that
$T_{a,b,n} \approx T_{a,b,-n}$. A Euclidean shear fixes the vertices
of the base and shows that every
$T_{a,b,n}$ is equivalent to some $T_{a',b',|n|}$, with $0 \le a',b'
\le |n|-1$. Given $(a,b,n)$, we define 
$$
c =1-a-b.
$$ 
Note that $c \equiv 1 \mod n$ if and only if 
 $n\ | \ a+b$. As we shall see, $c$ is an ``equal partner'' of $a$ and $b$
in $T_{a,b,n}$.

If $T \approx T'$, then  $vol(T) = vol(T')$. However, equal volumes do
not imply equivalence, even for clean tetrahedra.

\begin{lemma}(\cite[Thm. 5.6]{rz},\ \cite[pp. 144-145]{ha},\
    \cite[Thm. 5.1]{k})
We have $T_{a,b,n} \approx T_{d,e,n'}$ if and only if $|n| = |n'|$
and $d$ and $e$ are congruent $ \mod {|n|}$ to two of the elements in one
of the following triples:
\begin{equation}\label{E:equiv}
(a,b,c),\quad(a^{-1}, -ba^{-1}, -ca^{-1}),\quad (b^{-1}, -ab^{-1},
-cb^{-1}),\quad (c^{-1}, -ac^{-1}, -bc^{-1}).
\end{equation}
(If any of $\{a,b,c\}$ is not invertible mod $|n|$, then the
  corresponding triple does not appear in \eqref{E:equiv}.)
\end{lemma}
\begin{proof}
Since $f \in \mathcal L$ preserves volume, $|n| = |n'|$ is
necessary. Suppose $|n| = |n'|$, and after a possible  reflection,
suppose $n' = n > 0$. An affine map in $\rr^3$ is
determined by its values on the vertices non-degenerate tetrahedron;
since any affine map taking $T_{a,b,n}$ to $T_{d,e,n}$
will be volume-preserving, the issue is whether its coefficients are
integral. There are $4!=24$ cases. If the base of
$T_{a,b,n}$ is mapped to the base of $T_{d,e,n}$ by a
map in $\mathcal L$, then $f$ permutes $(0,0,0), (1,0,0), (0,1,0)$,
and so $f$ fixes $z$ and sends  $(x,y)$ to two of $\{x,y,1-x-y\}$.  
In this case, $(d,e)$ is congruent, mod $n$, to two
of the elements $(a,b,c)$, in some order.

We do one case to stand for the remaining 18: If $f(0,0,0) = (0,0,0)$,
$f(1,0,0) = (1,0,0)$, $f(0,1,0) = (d,e,n)$ and  $f(a,b,n) = (0,1,0)$,
then 
$$
f(x,y,z) = \left( x+dy - \left(\frac{a+bd}n\right) z,\  ey +
\left(\frac{1-eb}n \right)z,\ ny - bz \right).
$$
Observe that $f \in \mathcal L$ if and only if $\gcd(b,n) = 1$, $d
\equiv -ab^{-1}$ mod $n$ and 
$e \equiv b^{-1}$ mod $n$. The other cases are numbingly similar. 
\end{proof}
If we suspend $\rr^3$ in a hyperplane of $\rr^4$ via $(x,y,z) \mapsto
(1-x-y-z,x,y,z)$, then
$T_{a,b,n}$ is sent to the tetrahedron with vertices $(1,0,0,0),
(0,1,0,0), (0,0,1,0)$ and $(1-a-b-n,a,b,n)$. The six maps which
permute the base correspond to permutations of the first three
coordinates in $\rr^4$, and we can see directly that
$T_{a,b,n} \approx T_{a,c,n} \approx T_{b,a,n} \approx
T_{b,c,n} \approx T_{c,a,n} \approx T_{c,b,n}$.

\section{Clean and Empty Lattice Tetrahedra}

The following criterion was first studied by Reeve \cite{ree} in 1957.
 
\begin{lemma} Suppose $T$ is a non-degenerate lattice tetrahedron and
  suppose the only lattice points on the face $v_1v_2v_3$ are its
  vertices.  Then there exists
  $f \in \mathcal L$ so that $f(v_1) = (0,0,0)$,  $f(v_2) = (1,0,0)$
  and  $f(v_3) = (0,1,0)$. Thus, $T \approx T_{a,b,n}$ for some $(a,b,n)$
    with $n \ge 1$.
\end{lemma}
\begin{proof}
We construct the equivalence explicitly.
First translate by the first vertex, so that $v_1 =
(0,0,0)$; say that $v_2$ is now $(r,s,t)$. There are no lattice points
on the open edge $v_1v_2$, hence $1 = \gcd(r,s,t) =
\gcd(r,\gcd(s,t))$. Let $g = \gcd(s,t)$. Tweak the
last two coordinates, sending $v_2$ to $(r,g,0)$, and then tweak
the first two coordinates, so 
$(r,g,0) \mapsto (1,0,0)$. These tweaks fix $v_1$, and at this point,
we have $v_1 = (0,0,0)$ and $v_2 = (1,0,0)$.
 Suppose now that $v_3 = (i,j,k)$; again
$\gcd(i,j,k) = 1$; again tweak the last two coordinates so that
$(i,j,k) \mapsto (i,q,0)$, fixing $v_1$ and $v_2$. A Euclidean shear 
fixes $v_1, v_2$ and sends $v_3 = (i,q,0)$ to $(p,q,0)$, where $p
\equiv i \mod q$ and $0 \le p \le q-1$. Since $1 =
\gcd(p,q)$, if $q > 1$ then $p \ge 1$.
We claim that $q=1$, so $p=0$.  Suppose otherwise that $q \ge 2$.
Note that the non-vertex lattice point
$$
(1,1,0) = \left(\frac{p-1}q\right) \cdot (0,0,0) +  \left(\frac{q-p}q
\right)\cdot (1,0,0) +  \left(\frac 1q\right) \cdot (p,q,0)
$$
is on the face $v_1v_2v_3$, violating the
cleanliness of $T$. Therefore, $q=1$. (This last step can be replaced by
an appeal to Pick's Theorem; see \cite[p.233]{rz}.)
\end{proof}

\begin{theorem} (\cite[pp.389-390]{ree})
The lattice tetrahedron $T$ is clean if and only
  $T \approx T_{0,0,1}$ or 
$T \approx T_{a,b,n}$, where 
\begin{equation}\label{E:relpr}
n \ge 2,\quad  0 \le a,b \le n-1,\quad\gcd(a,n) = \gcd(b,n) = \gcd(c,n) = 1.
\end{equation}
 This equivalence can be effected with $f \in \mathcal L$ sending
 $v_1$, $v_2$, $v_3$, in that order, to $(0,0,0)$, $(1,0,0)$ and $(0,1,0)$.
\end{theorem}
\begin{proof}
First, observe that $T_{0,0,1}$ is clean. Suppose
$T=T_{a,b,n}$, where \eqref{E:relpr} holds. Then the face $v_1v_3v_4$
is contained in the plane $nx=az$, and if $w = (r,s,t)$ is a lattice 
point on this face, then $nr = at$. Since $\gcd(a,n) = 1$, it follows
that $n\ | \ t$. But  $v_1v_3v_4$ lies between the planes 
$z=0$ and $z=n$, so $0 \le t \le n$ and hence
$t = 0$ (so $r=0$ and $w = v_1$ or $v_3$) or $t=n$ (so $r = a$ and
$w=v_4$.) It follows that no non-vertex lattice points lie on the face
$v_1v_3v_4$. Since $v_1v_2v_4$ is contained in the plane $ny = bz$ and
$v_2v_3v_4$ lies in the plane $n(x+y-1) = (a+b-1)z$, similar
arguments, applied to $\gcd(b,n) = 1$ and $\gcd(c,n) = \gcd(1-a-b,n) =
1$ show that $T$ is clean.

Conversely, suppose $T$ is clean and $vol(T) = n/6$. Apply Lemma
2 so that $T \approx T_{a,b,n}$, with $0 \le a,b \le n-1$. If $n=1$,
then  $T \approx T_{0,0,1}$. Otherwise, $n \ge 2$. We need to show that
$\gcd(a,n) = \gcd(b,n) = \gcd(c,n) = 1$.  Suppose $g = \gcd(a,n) >
1$, and write  $(a,n) = (ga',gn')$. If $b = gb'$ for $b'
\in \zz$, then $v_4 = (ga',gb',gn')$ and the lattice point
$(a',b',n')$ is on the open edge $v_1v_4$, which 
is impossible. Accordingly, $b/g \notin \zz$;  let $m = \lfloor
b/g \rfloor$ and  $\ell = b - gm$, noting that $1 \le \ell \le g-1$.  
Observe that the non-vertex lattice point
$$
(a',m+1,n') = \left( \frac{\ell-1}g\right)\cdot(0,0,0)  +
\left(\frac{g-\ell}g\right)\cdot (0,1,0) +
 \left(\frac 1g\right)\cdot(ga',b,gn')
$$
is on the face $v_1v_3v_4$, a contradiction. Similar arguments apply
to $b$ and $c$.
\end{proof}

The number of interior lattice points in $T_{a,b,n}$ for $n \ge 2$
requires a non-trivial computation, in contrast to the $n$ points in the
corresponding fundamental parallelepiped. The following result can be
pieced together from Theorems 4.5, 4.7 and 5.2 in \cite{rz}; we prove
it here directly.
\begin{theorem}
(i)  Suppose $T = T_{a,b,n}$ is clean (so \eqref{E:relpr} holds.) Let 
\begin{equation}\label{E:frct}
A_t: = 
\left\{ \frac{t(n-c)}n \right\} + \left\{ \frac{t(n-a)}n \right\} +
\left\{ \frac{t(n-b)}n \right\} +  \left\{ \frac tn \right\}.
\end{equation}
Then
$i(T) = \# \{t: 1\le t \le n-1 \text{ and } A_t = 1\}$.

(ii) Suppose $T$ is a clean tetrahedron and $w \in int(T) \cap
\zz^3$. Then $BC(w) =
\left(d_1/N,d_2/N,d_3/N,d_4/N\right)$ for positive integers $d_j$,
$\sum d_j = N$, so that $\gcd(d_j,N) = 1$. 

(iii) Given positive integers $d_j$, $\sum d_j = N$, so
that $\gcd(d_j,N) = 1$, let $\la
=\left(d_1/N,d_2/N,d_3/N,d_4/N\right)$.
Then there exists a clean tetrahedron $T$ with at least one 
interior lattice point $w$ for which $BC_T(w) = \la$.
\end{theorem} 
\begin{proof}
  
We first remark that if $w = (r,s,t) \in \zz^3$ and $T = T_{a,b,n}$,
then a routine computation shows that 
\begin{equation}\label{E:bc}
BC_T(w) = \left( 1 - r - s + \frac{t(a+b-1)}n,\ r-\frac{ta}n,\
s-\frac{tb}n,\ \frac tn\right).
\end{equation}
Observe for later use that  $\{\frac{-ta}n\} =
\{\frac{t(n-a)}n\}$, $\{\frac{-tb}n\} = \{\frac{t(n-b)}n\}$ and 
$\{\frac{t(a+b-1)}n\} = \{\frac{t(n-c)}n\}$, because in each case, the
two fractions differ by an integer.

(i)  If $w \in int(T)$,  \eqref{E:bc}
implies that  $1 \le t \le n-1$,
and since $\la_2(w), \la_3(w)$ lie in $(0,1)$, we must have $r: = r(t) =
-\lfloor\frac{-ta}n\rfloor$ and $s := s(t) = -\lfloor\frac{-tb}n\rfloor$, so
that $\la_2(w) = \{\frac{t(n-a)}n\}$ and $\la_3(w) =  \{\frac{t(n-b)}n\}$.
Thus the only possible interior points have the shape $w_t:=
(r(t),s(t),t)$, and $w_t \in int(T)$ if and only if
$\la_1(w_t) \in (0,1)$. Since $\la_1(w_t) = 1 - \sum_{j=2}^4
\la_j(w_t)$ differs from $\frac{t(n-c)}n$ by an integer, we see
that $\la_1(w_t) \in (0,1)$ if and only if it actually equals
$\{\frac{t(n-c)}n\}$; that is, $w_t \in T$ if and only if $A_t = 1$.

(ii) The unimodular map taking $T$ to $T_{a,b,n}$ maps $w$ to
 $w_t$ for some $t$ with $1 \le t \le n-1$; let  $g = \gcd(t,n)$. By
taking reduced fractions in  
\begin{equation}\label{E:bcint}
BC(w_t) = \left( \left\{ \frac{t(n-c)}n \right\}, \left\{
\frac{t(n-a)}n \right\},
\left\{ \frac{t(n-b)}n \right\},  \left\{ \frac tn \right\} \right), 
\end{equation}
and recalling \eqref{E:relpr}, 
we see that \eqref{E:bcint} gives the desired shape with  $N = n/g$.

(iii) Since $\gcd(d_4,N) = 1$, we can choose  $m \in \zz$ so that
$md_4 = 1 + sN$. Now, let $T = T_{N-m d_2,N-md_3,N}$ and $w = (d_4 -
sd_2,d_4 - sd_3,d_4)$. It is easy to check that $BC_T(w) = \la$.
\end{proof}

There is an important clean tetrahedron in which the only interior point has
$g > 1$. Let $T=T_{3,7,20}$; since 3, 7 
and $c(3,7) = -9$ are relatively prime to 20, $T_{3,7,20}$ is clean.
A routine computation shows that
$\{ \frac{9t}{20}\} + \{ \frac{17t}{20} \} +
\{ \frac{13t}{20} \} +  \{ \frac t{20} \} = 1$ for $1 \le t \le 19$
if and only if $t=5$, so $T$ is a 1-point tetrahedron with $BC(w) = 
(\{ \frac{45}{20}\}, \{ \frac{85}{20} \},
\{ \frac{65}{20} \} ,  \{ \frac {5}{20} \}) = ( \frac 14,  \frac 14,
\frac 14,  \frac 14)$. We show in Theorem 14 that, up to equivalence,
$T$  is the only 1-point tetrahedron for which $g > 1$.

Let $\Lambda = \bigl\{ \sum_{j=1}^4 r_jv_j : r_j \in \zz,\  \sum_{j=1}^4
r_j = 1 \bigr\}$ 
denote the affine lattice determined by the vertices of
$T$. When $T = T_{a,b,n}$,
the fundamental region of $\zz^3/\Lambda$ can be taken to be a
parallelepiped with edges $v_j - v_i$, $v_k - v_i$, $v_\ell - v_i$,
for any permutation of the vertices. It follows from  \eqref{E:bc}
that $\{w_t : 0 \le t \le n-1\}$ is a set of representatives of
$\zz^3/\Lambda$. Let $M(a,b,n)$ denote the $(n-1)
\times 4$ matrix whose $t$-th row is $BC(w_t)$; namely,
$\left(\left\{ \tfrac{t(n-c)}n \right \}, \left\{
\tfrac{t(n-a)}n\right \}, \left\{  \tfrac{t(n-b)}n\right \}, \tfrac tn 
\right)$.  
Permutation of the
first three columns yields $M(a,c,n), M(c,b,n)$, etc. In order to
involve the fourth column, we must permute the rows as well to insure
that the last entry in the first row is $1/n$. If,  for example, $u \equiv
-b^{-1} \mod n$, then the $u$-th row of $M(a,b,n)$ is 
$$
\left( \tfrac{n-(-b^{-1}c)}n, \tfrac{n-(-b^{-1}a)}n, \tfrac 1n, 
\tfrac {n-b^{-1}}n \right),
$$
and after permuting the last two columns, we
obtain the first row of $M(-b^{-1}a, b^{-1},n)$. If $gcd(u,n) =
1$, then $\{ut \mod n: 1 \le t \le n-1\}$ is simply a permutation of 
 $\{t \mod n: 1 \le t \le n-1\}$, so replacing $t$ by $ut$ in the
definition of $M(a,b,n)$ permutes its rows. It follows in this way
from Lemma 1 that $T_{a,b,n} \equiv T_{a',b',n}$ if and only if
$M(a',b',n)$ can be derived from $M(a,b,n)$ after a permutation of rows
and columns. 

As an application of Theorem 4(i), consider $T_{n-1,n-1,n}$; since
$c=3-2n$, we assume   $\gcd(n,3) = 1$. In this case, $BC(w_t) =(
\{\frac{-3t}n\},\frac tn, \frac tn, \frac tn)$, and $A_t=1$ precisely
for $1 \le t \le \lfloor \frac n3 \rfloor$, so $i(T) = 
\lfloor \frac n3 \rfloor$. Han Duong has computed $i(T)$ for all clean
$T_{a,b,n}$ with $1 \le n \le 100$ and found that, for $k = i(T)$, the
inequality $3k+1 \le n \le 12k+8$ holds. In every case that $n = 3k+1$, the
tetrahedra are equivalent to $T_{3k,3k,3k+1}$; in every case that $n =
12k + 8$, the tetrahedra are equivalent to
$T_{2k+1,4k+3,12k+8}$. Additionally, he has shown that
$i(T_{2k+1,4k+3,12k+8}) = k$ for all $k \ge 1$. Duong and the author
conjecture that these bounds are, in fact, sharp for all $k \ge 1$.
Since $A_t + A_{n-t} = 4$, in the proof of Theorem 4, 
it follows that $n \ge 2k+1$ in any case. Heuristically, one would ``expect''
$k \approx \frac n6$ if the lattice points were evenly spaced in the
fundamental parallelepiped; the conjecture suggests bounds of roughly
$\frac n3$ and $\frac n{12}$. This conjecture will be discussed at greater
length in \cite{dhpr}.

Duong has an elegant proof of the bound $n=3k+1$, though without
uniqueness. Suppose $T$ is a lattice tetrahedron
(not necessarily clean) with $k$ interior points. We use these points
to subdivide $T$. Every point used in this subdivision is
either interior to one of the tetrahedra, or on an interior face or on an
interior edge. In these cases, the subdivision creates 3, 4 and 3 new
tetrahedra, respectively, and so in the end, $T$ is a union of at
least  $3k+1$ lattice tetrahedra, each of
which has volume $\ge 1/6$.

The characterization of empty tetrahedra was first made by White
\cite[pp.390-394]{w}, using a longish combinatorial
proof. P.~Noordzij \cite{n} gave a generalization with a longer,
elementary proof. There have since been a short, but sophisticated
proof using $L$-functions by D.~Morrison and G.~Stevens
\cite[pp. 16-17]{ms},  and
combinatorial proofs by H.~E.~Scarf \cite[pp. 411-413]{s} (based in
part on work of
R.~Howe) and Handelman \cite[pp. 145-148]{ha}. We present yet another 
 elementary proof. 

\begin{theorem}
 The clean lattice tetrahedron $T$ is empty
if and only if $T \approx T_{0,0,1}$  or $T \approx T_{1,d,n}$,
where $\gcd(d,n) = 1$ and $1 \le d \le n-1$. This equivalence can be
effected by a unimodular map sending $v_4$ to $(1,d,n)$.
\end{theorem}
\begin{proof}
If $T$ is empty and $vol(T) = \frac 16$, then $T \approx T_{0,0,1}$,
which is clearly empty.

We may now  assume $vol(T) = \frac n6 > \frac 16$, so $T \approx
T_{a,b,n}$, $n \ge 2$.
Observe that  $T_{1,d,n}$ is contained in the
slab $0 \le x \le 1$. Thus, if $w = (r,s,t) \in T_{1,d,n} \cap \zz^3$,
then $r = 0$ (and $w$ is on 
the edge $v_1v_3$) or $r=1$ (and $w$ is on
the edge from $v_2v_4$.) Since $\gcd(d,n) = 1$, this
implies that $w$ is a vertex, and so $T_{1,d,n}$ is empty. (Put another
way, $T_{1,d,n}$ has lattice width 1, because it lies in the
consecutive planes $x=0$ and $x=1$; White expressed his result by
saying that the vertices of $T$ lie in consecutive lattice planes. We
return to this topic in section five.)

Conversely, suppose $T_{a,b,n}$ is empty, where $n \ge 2$. Referring
to \eqref{E:frct}, it follows from Theorem 4(i) that $A_t > 1$ for
$1 \le t \le n-1$. Further,
$$
A_t + A_{n-t} = 
\left( \left\{ \frac{t(n-c)}n \right\} + \left\{ \frac{(n-t)(n-c)}n
\right\} \right)  + \cdots
$$
can be arranged into a sum of four pairs of terms each of which sums to 1,
hence $A_t + A_{n-t} = 4$, and so $A_t = 2$ for $1\le t \le n-1$. It
is now convenient to let $n-1 \ge a_3 \ge a_2 \ge a_1$ be the ordered
rearrangement of $\{n-a,n-b,n-c\}$. In view of Lemma 1, it suffices to 
show that  $a_3 = n-1$.
Writing $a_0 = 1$, we have $A_t = \sum_{j=0}^3  \{ ta_j/n\} = 2$
for $1 \le t \le n-1$, and taking $t=1$, we see that  
 $a_0+a_1+a_2+a_3 = 2n$.
Define $B_t: = \sum_{j=0}^3  \lfloor ta_j/n \rfloor$.
Then 
$$
A_t + B_t =\sum_{j=0}^3\left( \left\{ \frac{ta_j}n \right\}  + \left\lfloor
\frac{ta_j}n\right\rfloor\right)   = \sum_{j=0}^3 \frac{ta_j}n = t(A_1
+ B_1),  
$$
and since $A_1 = 2$ and $B_1 = 0$, we conclude that $B_t = 2(t-1)$ for
$1 \le t \le n-1$. In particular, $B_2 = 2$, and since $\lfloor
2a_j/n \rfloor \in (0,2)$, we must have $a_3 \ge a_2 \ge n/2 >
a_1$. Write $a_3 = n - b_2$ and $a_2 = n - b_1$, so that  $b_2 \le b_1 <
n/2$ and note that $(n-b_1)+(n-b_2) + a_1 + 1 = 2n$, hence $b_1 + b_2
= a_1 + 1$.

We need to show that $b_2 = n-a_3 = 1$. Suppose not; then $2 \le b_2
\le b_1 \le a_1 - 1 < n/2$. 
We have
\begin{equation}\label{E:mess}
2(t-1) =  \left\lfloor  \frac{t(n-b_1)}n \right\rfloor +  \left\lfloor
\frac{t(n-b_2)}n \right\rfloor +  \left\lfloor  \frac{ta_1}n
\right\rfloor + \left\lfloor \frac tn \right\rfloor.
\end{equation}
But $t/n \in (0,1)$ and  $t(n-b_j)/n \notin \zz$ implies that
 $\lfloor  t(n-b_j)/n \rfloor = t +
 \lfloor  t(-b_j)/n \rfloor = t -  \lfloor
 tb_j/n \rfloor - 1$. Thus,
it follows from \eqref{E:mess} that
$$
2(t-1) = (t-1) -  \left\lfloor \frac{tb_1}n \right\rfloor +  (t-1) -
\left\lfloor \frac{tb_2}n \right\rfloor +  \left\lfloor
\frac{ta_1}n\right\rfloor,
$$
and so
\begin{equation}\label{E:key}
C_t := \left\lfloor \frac{ta_1}n\right\rfloor =  \left\lfloor
\frac{tb_1}n \right\rfloor +  \left\lfloor \frac{tb_2}n \right\rfloor.
\end{equation}

For $r \in \rr \setminus \zz$ and $k \in \zz$, let $\Del_k(r):=
\lfloor (k+1)r \rfloor 
- \lfloor kr \rfloor$. Since $\Del_k(r) = r + \{kr\} - \{(k+1)r\}$, we have
 $|\Del_k(r) - r| < 1$, hence  $\Del_k(r) =
\lfloor r \rfloor$ or $ \lfloor r \rfloor + 1$. In particular, if
$r \in (0,1)$, then  $\Del_k(r) \in \{0,1\}$.
Further,  $\Del_k(r)=1$ if and only if
there is an integer $h$ so that $kr < h \le (k+1)r$; that is, $k <
h/r \le k+1$, or $k+1 = \lceil h/r \rceil$. If we also know
that $h/r \notin \zz$, then  $k =  \lfloor
h/r \rfloor$, and so for fixed $r \in (0,1)$, $\Del_k(r) = 1$ on a sequence
of $k$'s with jumps of size
$\lfloor 1/r \rfloor$ or $\lfloor 1/r \rfloor + 1$.

For $1 \le t \le n-2$, let $D_t = C_{t+1} - C_t$. 
It follows from \eqref{E:key} that 
\begin{equation}\label{E:zap}
D_t = \Del_t\biggl(  \frac{a_1}n \biggr) = \Del_t\biggl(
\frac{b_1}n\biggr)  + \Del_t\biggl(\frac{b_2}n \biggr),
\end{equation}
so in particular, $D_t = 0$ or $1$.
Let $S = \{t : D_t = 1\}$. It follows from \eqref{E:zap} that
$$
S = \left\{\left \lfloor \frac n{a_1} \right\rfloor,\  \dots,\ 
\left\lfloor\frac {(a_1-1)n}{a_1} \right\rfloor\right \}. 
$$
(Observe that $kn/a_1 \not\in \zz$ for $1 \le k \le a_1-1$.) 
Letting $g = \lfloor n/a_1 \rfloor \ge 2$, we see that $S$ contains
$a_1-1$ integers, with jumps of $g$ or $g+1$. Let
$$
S_j =\left \{\left \lfloor \frac n{b_j} \right\rfloor,\  \dots,\ 
\left\lfloor\frac {(b_j-1)n}{b_j} \right\rfloor\right \}. 
$$
for $j=1,2$; it also follows from \eqref{E:zap} that   $S = S_1 \cup
S_2$, where $|S_j| = b_j - 1$; note that $a_1-1 = (b_1-1)+(b_2-1)$ and
$kn/b_j \not\in \zz$ for $1 \le k \le b_j-1$.
 Since $b_1 \ge b_2$, 
$\lfloor n/b_1 \rfloor \le \lfloor  n/b_2 \rfloor$ and so
$g =  \lfloor n/b_1 \rfloor$. Thus $S_1$ is also a sequence 
with jumps of size $g$ or $g+1$. Since $2g > g+1$, any elements in $S$
that are not in $S_1$ must come from the ends of $S$, not the
middle. But $S$ and $S_1$ have the same first element $g$ and the last
element, 
$$
\left\lfloor\frac {(a_1-1)n}{a_1} \right\rfloor =\left \lfloor n - \frac
n{a_1}\right \rfloor = n - (g + 1) = \left\lfloor\frac {(b_1-1)n}{b_1}
\right\rfloor, 
$$
which gives the contradiction. We conclude that $b_2 = 1$, completing
the proof. 
\end{proof}

Reeve \cite{ree} observed that $T_{1,d,n}$ is always empty, but
that $a=1$ or $b=1$ is not a necessary condition, as 
$T_{2,5,7}$ is empty.  In our notation, $c(2,5) \equiv 1 - 2 - 5
  \mod 7 = 1$. M.~Khan \cite{k} has given a formula for the number of
  equivalence classes of empty tetrahedra of volume $n/6$.

We conclude this section with an observation that will be
essential in the next section.

\begin{corollary} Suppose $T$ is a clean tetrahedron. Then $T$
  is empty if and only if, for every $w = (r,s,t) \in \zz^3$, the 
  $\la_{j,T}(w)$'s sum pairwise to integers. 
\end{corollary}
\begin{proof}
If $T$ is empty, then $T \approx T_{0,0,1}$ or $T \approx
T_{1,d,n}$ by Theorem 5. (Permutation of the vertices is
irrelevant to the condition given.)  For
$T = T_{1,d,n}$, and $w = (r,s,t)$,  \eqref{E:bc} shows that $\la_{1,T}(w) +
\la_{3,T}(w) = 1-r$ and $\la_{2,T}(w) + \la_{4,T}(w) = r$.

If $T$ is not empty and  $w$ is an interior point, then the relations 
$0 < \la_{j,T}(w) < 1$ and $\sum_j \la_{j,T}(w) = 1$ show that no sum of
two can be integral.
\end{proof}
The geometric interpretation of this result (see  \cite{k}) is
that in the fundamental region of $\zz^3/\Lambda$, every point in $\zz^3$
lies in one of the three interior diagonal parallelograms which avoid
$T$. Khan points out the surprising fact that there is no purely
geometric proof of this result. 
Alex Kasprzyk has pointed out to the author that clean tetrahedra
correspond to a cone which, in toric geometry, is a terminal quotient
singularity, and remarks that this situation is discussed in
\cite[p.379]{rei}. 

\section{1-point lattice tetrahedra}

The goal of this section is to prove the following classification
theorem, which has been proved in a somewhat different way by A. Kasprzyk
\cite{ka}. 

\begin{theorem} If $T$ is a 1-point lattice tetrahedron,
  then $T$ is equivalent to 
  $T_{3,3,4}$, $T_{2,2,5}$, $T_{2,4,7}$,
 $T_{2,6,11}$, $T_{2,7,13}$, $T_{2,9,17}$, $T_{2,13,19}$ or $T_{3,7,20}$.
\end{theorem}

The proof of this result will come from combining Theorems 13 and 14  below.
We begin our discussion with a derivation of the possible barycentric
coordinates for the interior lattice point. This was done in \cite{rz}
in a less transparent way. Our proof relies on two simple observations
about a clean 1-point lattice tetrahedron $T$ with interior lattice
point $w$. The first is that $BC_T(w)$ must have special arithmetic
properties (via Corollary 6). The second is that $w$ subdivides $T$
into four {\it empty} tetrahedra.  

First,
suppose $T = T(v_1,v_2,v_3,v_4)$ is a clean tetrahedron and $w
\in int(T)$, where
\begin{equation}\label{E:bas}
BC_T(w) = \la = (\la_1,\la_2,\la_3,\la_4) = \left( \frac{d_1}N,
\frac{d_2}N,  \frac{d_3}N, \frac{d_4}N \right).
\end{equation}
Since $T_{0,0,1}$ is empty, Theorem 3 implies that $T$ is equivalent
to some $T_{a,b,n}$, satisfying \eqref{E:relpr}, and we may assume
 $T = T_{a,b,n}$, where, by Theorem 4(ii), $n = gN$. (Caveat: in the
proof, we might subsequently apply
one of the maps from Lemma 1 to permute the vertices and replace
$T_{a,b,n}$ with an equivalent $T_{a',b',n}$. This, of course,
permutes the $d_j$'s.) 

The following is a slight restatement of Theorem 4.

\begin{theorem} 
(i) If $T$ is a 1-point lattice tetrahedron with interior point $w$
satisfying \eqref{E:bas}, then
\begin{equation}\label{E:nogen}
2 \le s \le N-1 \implies \sum_{j=1}^4 \left\{ \frac{sd_j}N  \right\} > 1.
\end{equation}
(ii) If $g=1$ and \eqref{E:nogen} holds, then $T$ is a 1-point lattice
tetrahedron. 
\end{theorem}
\begin{proof}
(i) Suppose $T$ is a 1-point lattice tetrahedron, but \eqref{E:nogen}
  fails for $s= s_0 \ge 2$. Then  
$$
w':= \sum_{j=1}^4 \{ s_0\la_j\}v_j = \sum_{j=1}^4 ( s_0\la_j - \lfloor
s_0\la_j 
 \rfloor) v_j =
s_0w - \sum_{j=1}^4 \lfloor s_0\la_j \rfloor v_j \in \zz^3
$$
is also an interior lattice point in $T$, hence $w' = w$, and 
by the uniqueness of barycentric coordinates, $\{s_0\la_j\} = \la_j$.
It follows that  $(s_0-1)\la_j \in \zz$ and $N \ | \ (s_0-1)d_j$, so
$N$ divides $s_0-1$, a contradiction.

(ii) Let $w = (r,s,t_0)$. Then by \eqref{E:bcint} and \eqref{E:bas},
\begin{equation}\label{E:link}
\left( \frac{d_1}N, \frac{d_2}N,  \frac{d_3}N, \frac{d_4}N \right)
=  \left( \left\{ \frac{t_0(n-c)}n \right\}, \left\{
\frac{t_0(n-a)}n \right\},
\left\{ \frac{t_0(n-b)}n \right\},  \left\{ \frac {t_0}n \right\} \right). 
\end{equation}  
As noted earlier, since $\gcd(t_0,N) = 1$, multiplication by $t_0$
permutes the non-zero 
residue classes mod $N$.
Thus, \eqref{E:nogen} and \eqref{E:link} imply that $A_t = 1$ only for
$t=t_0$, hence $k=1$ and $T$ is a 1-point lattice tetrahedron.  
\end{proof}

The condition $g=1$ is essential in Theorem 8(ii). For example,
consider $T_{7,7,8}$ and $w = (2,2,2)$. Then $BC(w) = (\frac 14,
\frac 14, \frac 14, \frac 14)$, which clearly satisfies
\eqref{E:nogen} with $g=2$. However, $(1,1,1) = \frac 12(v_1 + w)$ is another
interior point in $T_{7,7,8}$.

The second observation we make is simpler. The interior
point $w$ 
subdivides $T$ into four tetrahedra: $T_1 = T(v_2,v_3,v_4,w)$, $T_2 =
T(v_1,v_3,v_4,w)$, $T_3 = T(v_1,v_2,v_4,w)$ and $T_4 = T(v_1,v_2,v_3,w)$.
(We fix these notations for the rest of this section.) The
following lemma was used in \cite{m}, though without the simplified
implications of Lemma 6.

\begin{lemma}
The four lattice tetrahedra $T_1$, $T_2$, $T_3$ and $T_4$ 
are empty.
\end{lemma}
\begin{proof}
This is immediate, because $T_j \cap \zz^3 \subseteq T \cap \zz^3 =
\{v_1, v_2, v_3, v_4, w\}$. 
\end{proof}

\begin{lemma}
Suppose $T$ is a 1-point tetrahedron, then in
the notation of \eqref{E:bas}, each $d_i$ divides one of  $d_j+d_k$, $d_j +
d_\ell$ or $d_k+d_\ell$, where $\{i,j,k,\ell\} = \{1,2,3,4\}$.
\end{lemma}
\begin{proof}
Suppose $i = 4$ for clarity. Then since $w = \sum_{j=1}^4 \la_j
v_j$, we have 
\begin{equation}\label{E:flip}
v_4 = 
\left(\frac{-d_1}{d_4}\right)v_1 + \left(\frac{-d_2}{d_4}\right)v_2
+ \left(\frac{-d_3}{d_4}\right)v_3 + \left(\frac{N}{d_4}\right)w,
\end{equation}
so $BC_{T_4}(v_4) =
(-d_1/d_4,-d_2/d_4,-d_3/d_4,N/d_4)$. By Lemma 9, $T_4$ is
empty, and hence by Corollary 6, $d_4$ divides one of
$\{d_1+d_2,d_1+d_3,d_2+d_3\}$.
\end{proof}
\smallskip
We shall say that  $(d_1,d_2,d_3,d_4) \in \nn^4$ is \emph{ripe} if,
for  $N = \sum d_i$, we have:
\begin{enumerate}
\item[(i)] For each $j$, $\gcd(d_j,N) = 1$.

\item[(ii)] Each $d_i$ divides one of  $d_j+d_k$, $d_j +
d_\ell$ or $d_k+d_\ell$.

\item[(iii)] If $d_j= d_i$ or $d_j = 2d_i$, then $d_i
  = 1$.  

\end{enumerate}
\begin{lemma}
If $T$ is a 1-point lattice tetrahedron, then $(d_1,d_2,d_3,d_4)$
is ripe. 
\end{lemma}
\begin{proof}
Conditions (i) and (ii) follow from Theorem 4(ii) and Lemma 10. For
(iii), suppose for concreteness that $d_1 > 1$ and $d_2 = \ep d_1$, where 
$\ep \in \{ 1,2\}$. 
Let $1 \neq s \equiv d_1^{-1} \mod N$. We claim that for $j = 3,4$,
$sd_j \not\equiv N-1 \mod N$. Suppose otherwise. Then $sd_j \equiv N-1
\mod N$ implies  
$d_j \equiv N-d_1 \mod N$, so $d_j = N-d_1$, hence $d_1+d_2+d_3+d_4 > N$.
It follows that  $\{sd_j/N\} \le (N-2)/N$ for $j=3,4$, thus
$$
\sum_{j=1}^4 \left\{ \frac{sd_j}N \right \} = \frac 1N + \frac \ep N +
\sum_{j=3}^4 \left\{ \frac{sd_j}N \right \} \le \frac{1+\ep +2(N-2)}N
< 2. 
$$
Since this sum is an integer, it must equal 1, which contradicts
Theorem 8(i).
\end{proof}
We now need a tedious case-analysis, which is nevertheless shorter than the
tedious case-analysis in the corresponding proof in \cite{rz}.
\begin{theorem}
If $(d_1,d_2,d_3,d_4)$ is ripe and $d_1\le d_2\le d_3\le d_4$, then
$(d_1,d_2,d_3,d_4)$ is $(1,1,1,1)$, $(1,1,1,2)$,
$(1,1,2,3)$, $(1,2,3,5)$, $(1,3,4,5)$, $(2,3,5,7)$ or $(3,4,5,7)$.
\end{theorem}
\begin{proof}
It is easy to verify that each of the quadruples given in the
statement is ripe.
Note also that if $d_j = \al_j t + \be_j u$ for integers $\al_j,
\be_j, t, u$, $1 \le j \le 4$,
 then $\gcd(t,u) = 1$; otherwise we would have $\gcd(d_j,N) > 1$,
violating (i).

We first consider the cases in which at least two $d_j$'s are
equal, so  $d_1 = d_2 = 1$ by (iii).
 If $(1,1,1,u)$ is ripe, then $u \ | \ 2$ by (ii), so $u = 1$ or 2. These
 give the first two examples. If
$(1,1,t,u)$ is ripe, with $t,u \ge 2$, then $u$ divides 2 or $t+1$ 
by (ii). But $u = 2$ implies $t=2$ by monotonicity,
 violating (iii). Therefore, $u \ | \
t+1$, and $u \ge t$ implies $u = t+1$. But if $(1,1,t,t+1)$ is ripe,
then $t>1$ divides 2 or $t+2$, so $t=2$, giving $(1,1,2,3)$,
the third example. 

We now assume 
\begin{equation}\label{E:monot}
1 \le d_1 < d_2 < d_3 < d_4.
\end{equation}
 Since $d_4$ divides a pair-sum less than $2d_4$, it must equal that
 pair-sum. There are three cases, which we consider in turn:
 (a) $d_4 = d_2 + d_3$, (b) $d_4 = d_1 +
d_3$, (c) $d_4 = d_1 + d_2$.

If  $(d_1,d_2,d_3,d_2+d_3)$ is ripe, then $d_3$ divides
 $d_1+d_2,d_1+d_2+d_3$ or $2d_2+d_3$ and hence $d_1+d_2$ or $2d_2$.
Both  are $< 2d_3$, so either $d_3 = d_1 + d_2$ or $d_3 = 2d_2$. The
latter implies that $d_2 = 1$ by (iii), violating \eqref{E:monot}, so
$d_3 = d_1 + d_2$.

After writing $d_1 = t$ and $d_2 = u$, we now
suppose that $(t,u,t+u,t+2u)$ is ripe, where $u \ge 2$.  Then  $u$ divides 
$2t+u$, $2t+2u$ or $2t+3u$, and so $u \ | \ 2t$. Since $\gcd(t,u) = 1$,  
$u = 2$ and  $(t,2,t+2,t+4)$ is ripe. It follows from \eqref{E:monot}
that  $2=t+1$ by \eqref{E:monot}, and we obtain $(1,2,3,5)$, the
fourth example.  This completes case (a).
 
Next, suppose that   $(d_1,d_2,d_3,d_1+d_3)$ is ripe, so $d_3$ divides
 $d_1+d_2$, $2d_1+d_3$ or $d_1+d_2+d_3$, and hence either $d_1+d_2$ or
$2d_1$. Both are $< 2d_3$, so either $d_3 = d_1 + d_2$ or $d_3
= 2d_1$. Again,  $d_3 = 2d_1$ implies $d_1 = 1$ by (iii), so $d_3 = 2$,
violating \eqref{E:monot}.  
Thus, $d_3 = d_1 + d_2$ and $(t,u,t+u,2t+u)$ is ripe, so again,
$u \ge 2$. 
 We have $N=4t+3u$, so (i) implies that $u$ is odd. 
Since $t$ divides  $t+2u$, $2t+2u$ or $3t+2u$, we
have $t \ | \ 2u$, so $t \in\{ 1,2\}$. If $t=1$ and  $(1,u,1+u,2+u)$ is
ripe, then  $u$  divides 
$2+u$, $3+u$ or $3+2u$, forcing $u=3$ and $(1,3,4,5)$, the fifth example.
 If $t=2$ and $(2,u,2+u,4+u)$ is ripe, then $u$ divides
$4+u$, $6+u$ or $6+2u$. Again, $u=3$ is the only choice, giving
$(2,3,5,7)$, the sixth example.
 This completes case (b).

Finally, suppose $(d_1,d_2,d_3,d_1+d_2)$ is ripe, so $d_3$ divides
$d_1+d_2$, $2d_1+d_2$ or $d_1+2d_2$. Since $d_1 + d_2 < 2d_3$,
these are each $< 3d_3$.
If $d_3 =d_1+d_2$, $2d_1+d_2$ or $d_1+2d_2$, then 
 $d_3 \ge d_4$, violating \eqref{E:monot}.
Thus, $2d_3$ equals $d_1+d_2, 2d_1 + d_2$ or $d_1+2d_2$, and the first
is impossible, so either  $2d_3 = 2d_1+d_2$ or $2d_3 = d_1+2d_2$, so
 either $d_2$ or $d_1$ is even and either $(t,2v,t+v,t+2v)$ or
$(2v,t,t+v,t+2v)$ is ripe. Except for the order of the first two
terms, these are the same case.
Since $2v$ divides $2t+v$, $2t+2v$ or $2t+3v$, 
 $v \ | \ 2t$, so $v \in \{1,2\}$. Suppose $v=1$ and $(t,2,t+1,t+2)$ or
$(2,t,t+1,t+2)$ is ripe. The former is impossible by \eqref{E:monot};
the latter implies that $t \ge 3$ divides $t+3$, $t+4$ or $2t+3$,
which implies $t=3$ (and $t+1 = 4$) or $t=4$. In either case, 
(iii) is violated. In the final case, $v=2$ and $(t,4,t+2,t+4)$ or
$(4,t,t+2,t+4)$ is ripe. The former implies that $4 = t+1$, yielding
$(3,4,5,7)$, the final example. The latter implies that $t > 4$ is odd
and divides $t+6$, $t+8$ or $2t+6$, which is impossible. 
 This completes case (c) and the proof.
\end{proof}

The concept of ripeness seems worth considering in its own right,
especially if 
(iii) is jettisoned. Without (iii), there are several additional
families of  ripe quadruples satisfying (i) and (ii):
 $(1,t,t,t)$, $(2,t,t,t,)$ (where $\gcd(t,2) = 1$),
$(1,t,t,2t)$, $(3,t,t,2t)$ (where $\gcd(t,3) = 1$), $(1,t,2t,3t)$ and
$(5,t,2t,3t)$ (where $\gcd(t,5) = 1$). These families generalize the first
four of the examples given. It would also be interesting to study
ripeness in $n$-tuples for $n\ge 5$.

We now use ripeness to give a new proof of Theorem 5.9 in
\cite{rz}. 

\begin{theorem}
(i) If  $T$ is a 1-point lattice tetrahedron
with interior point $w$, then, up to a permutation of the coordinates, 
$BC_T(w)$ is one of the following quadruples:
 $\la^{(1)} = \left( \frac 14,  \frac 14,\frac 14,
\frac 14  \right)$,  $\la^{(2)} = \left( \frac 15,  \frac 15,\frac 15,
\frac 25  \right),$ $\la^{(3)} = \left( \frac 17,  \frac 17,\frac 27,
\frac 37  \right),$ $\la^{(4)} = \left( \frac 2{11},  \frac
1{11},\frac 3{11}, \frac 5{11}  \right),$ $\la^{(5)} = \left( \frac
1{13},  \frac 3{13},\frac 4{13}, \frac 5{13} \right),$ $\la^{(6)} =
\left( \frac 2{17},  \frac 3{17},\frac 5{17}, \frac 7{17} \right),$
$\la^{(7)} = \left( \frac 3{19},  \frac
  5{19},\frac 4{19},  \frac 7{19} \right).$

(ii) If  $T$ is a 1-point lattice tetrahedron and $g=1$, then $T$ is
  equivalent to $T_{3,3,4}$, $T_{2,2,5}$, $T_{2,4,7}$,
 $T_{2,6,11}$, $T_{2,7,13}$, $T_{2,9,17}$ or $T_{2,13,19}$.  
\end{theorem}
\begin{proof}
(i) By Lemma 10, $BC_T(w)$ must be ripe. It is easily checked
  that every ripe quadruple in Theorem 12 satisfies \eqref{E:nogen},
  and so is $BC_T(w)$ for some 1-point tetrahedron by Theorem 8(ii).

(ii) Suppose $g=1$. Then $n=N$ and
$\la = \la^{(\ell)}$ up to some permutation
  of the coordinates. Take $f \in \mathcal L$ (if necessary) to permute the
vertices of $T$ so that $\la = \la^{(\ell)}$ and $T$ is a clean $T_{a',b',n}
(\approx T_{a,b,n})$ via Lemma 1. Then $a',b' \in \zz$ and
$$
w = \left( \frac{d_1}n\right) \cdot (0,0,0) +\left(
\frac{d_2}n\right)\cdot(1,0,0) + 
\left(\frac{d_3}n\right)\cdot(0,1,0) +
\left(\frac{d_4}n\right)\cdot(a',b',n). 
$$
Note that $w \in \zz^3$ if and only if 
 $d_2 + d_4a' \equiv d_3 + d_4b' \equiv 0 \mod n$; that is, $a' \equiv
-d_2d_4^{-1} \mod n$ and $b' \equiv -d_3d_4^{-1} \mod n$. Using
$(d_2,d_3,d_4)$ from the seven cases in (i), we obtain  the seven
tetrahedra in (ii). 
\end{proof}

An application of Lemma 1 yields the equivalent clean 1-point tetrahedra
of the form $T_{a,b,n}$ with $0 \le a,b \le n-1$. (In case
parameters are repeated, treat the triples as multi-sets.)
\smallskip
\begin{enumerate}

\item  $T_{3,3,4}$;

\item  $T_{a,b,5}$, $\{a,b\} \in \{2,2,2\}, \{3,4,4\}$;

\item $T_{a,b,7}$, $\{a,b\} \in \{2,2,4\}, \{2,3,3\}, \{4,5,6\}$;

\item $T_{a,b,11}$, $\{a,b\} \in \{2,3,7\}, \{2,4,6\}, \{3,4,5\},
  \{6,8,9\}$; 

\item $T_{a,b,13}$, $\{a,b\} \in \{2,3,9\}, \{2,5,7\}, \{3,4,7\},
  \{8,9,10\}$; 

\item $T_{a,b,17}$, $\{a,b\} \in \{2,3,13\}, \{2,7,9\}, \{5,6,7\},
  \{4,5,9\}$; 

\item $T_{a,b,19}$, $\{a,b\} \in \{2,5,13\}, \{3,4,13\}, \{3,7,10\},
\{4,5,11\}$.
\end{enumerate}
\smallskip
We complete our characterization of 1-point lattice tetrahedra.
Suppose now that  $T$ is a 1-point lattice tetrahedron and $g = n/N>
1$, where $BC_T(w) = \la^{(\ell)}$ for some $1 \le \ell \le 7$ and, without
loss of generality,
$\la^{(\ell)}_4$ is the largest component of $\la^{(\ell)}$.
 Since $vol(T_4) = \la_4 \cdot vol(T) = \la_4 n/6 = gd_4/6 > 1/6$, the
 (empty) tetrahedron $T_4$ is equivalent to $T_{1,e,gd_4}$, where 
\begin{equation}\label{E:cond1}
1  \le e \le  gd_4 - 1, \qquad \gcd(e,gd_4) = 1.
\end{equation}
 We apply a unimodular map
 so that 
\begin{equation}\label{E:base}
v_1=(0,0,0),\ v_2 = (1,0,0),\ v_3 = (0,1,0),\ w = (1, e, gd_4).
\end{equation}
In doing so, the vertices $v_1, v_2, v_3$ may be permuted, and so 
$BC_T(w) = \la^{(\ell)}$, with some uncertainty about the order of the
first three cooordinates.
 It follows from \eqref{E:flip} that
\begin{equation}\label{E:v}
v_4 =  \left( \frac{N-d_2}{d_4}, \frac{eN-d_3}{d_4}, gN \right).
\end{equation}
Since $v_4 \in \zz^3$, we see that $d_4$ must divide $N-d_2$ and
$eN-d_3$. This first condition reduces to $d_4 \ | \ d_1 + d_3$, which
has already been accommodated in Theorem 13(i); in fact, it is easy to
check that the first coordinate of $v_4$ is 3 (if $\ell = 1$), and 2
(if $\ell \ge 2$). The second condition is that $eN \equiv d_3$ mod
$d_4$. Further, since $T$ is clean, Theorem 3 implies that 
\begin{equation}\label{E:cond2}
\gcd\left(\frac{N-d_2}{d_4},gN \right) = \gcd\left(\frac{eN-d_3}{d_4},gN
\right) = \gcd\left(\frac{eN+d_1}{d_4},gN \right) = 1.
\end{equation}
(We have used here that the numerator of $-c = a+b-1$ is $N-d_2 + eN -
d_3 -d_4$.)

We discuss the seven cases in turn.

First, suppose $\ell = 1$, so that
$d_1=d_2=d_3=d_4=1, N=4$, and
$$
v_4 = (3,4e-1,4g), \qquad w = (1,e,g),
$$
where $1 \le e \le g-1$. Observe that the unimodular map $(x,y,z) \mapsto
(x, 1-x-y+z,z)$ permutes $v_1$ and $v_3$, fixes $v_2$
and sends $v_4 \mapsto (3,4(g-e)-1,4g)$ and $w \mapsto (1,g-e,g)$,
hence we may assume without loss of generality that $e \le
g/2$. Further, if $g$ is even, then $\frac 12(v_2 + w)$ or $\frac
12(v_4+w)$ is a lattice point, depending on whether $e$ is even or
odd; thus, $g$ is odd and $1\le e \le (g-1)/2$. Here, \eqref{E:cond1}
and \eqref{E:cond2} imply that $\gcd(e,g) = \gcd(3,4g) = \gcd(4e-1,g)
= \gcd(4e+1,4g) = 1$. In particular, $g \neq 3$.

If $g=5$, then the possible values for $e$ are 1 or 2, and the first
is ruled out by $gcd(5,20)>1$, hence $e=2$. In this case $T =
T_{3,7,20}$, which has 
already been identified as another 1-point tetrahedron. Now suppose
$g \ge 7$. We compute $BC_{T_j}(0,0,1)$ for $j = 2,3$: 
$$
(0,0,1) = 
\left( \frac{e+g+1}g\right) \cdot v_1 + \left(
\frac{-e+1}g\right) \cdot v_3 + \left( \frac 
1g\right) \cdot v_4 + \left( \frac{-3}g\right) \cdot w;
$$
$$
(0,0,1) =
 \left( \frac{2e+g}g\right) \cdot v_1 + \left(
\frac{e-1}g\right) \cdot v_2 +\left(  \frac 
eg\right) \cdot v_4 +\left(  \frac{-4e+1}g\right) \cdot w.
$$

Since $T_2$ and $T_3$ are empty, 
Corollary 6 implies that $g$ divides 2, $e-2$ or $e+2$ and $g$ divides
$2e-1$, $3e-1$ or $3e$. Since $e+2 \le \frac{g+3}2 < g$, the
first condition implies that $e=2$, and this means that the second
condition is impossible, and there are no ``new'' 1-point tetrahedra
with $\la = \la^{(1)}$. This reproduces the result of Mazur \cite{m}.

We now consider $\ell \ge 2$. In this case, we have $\frac{N-d_2}{d_4}
= 2$, and so by \eqref{E:cond2}, $\gcd(2,gN)=1$; thus, $g \ge 3$ is odd.

Suppose $\ell = 2$, so $d_1=d_2=d_3=1,d_4=2,N=5$. It follows from
\eqref{E:v} that $5e-1$ is even, hence $e = 2k+1$ is odd. Then
$$
w = (1,2k+1,2g), \qquad v_4 = (2,5k+2,5g),
$$
with $0 \le k \le g-1$ and  $\gcd(2k+1,2g) = \gcd(5k+2,5g) =
\gcd(5k+3,5g) = 1$. If $g=3$, note that $k=1,2,0$ (in that order) are
ruled out by the gcd conditions, so
we must have $g \ge 5$. Again, we compute 
 $BC_{T_j}(0,0,1)$ for $j = 2,3$:
$$
(0,0,1) = \left( \frac{g+k+1}g\right) \cdot v_1 + \left(
\frac{-k}g\right) \cdot v_3 + \left( \frac 
1g\right) \cdot v_4 + \left( \frac{-2}g\right) \cdot w;
$$
$$
(0,0,1) =  \left( \frac{g+2k+1}g\right) \cdot v_1 + \left(
\frac{k}g\right) \cdot v_2 +\left(  \frac 
{2k+1}g\right) \cdot v_4 +\left(  \frac{-5k-2}g\right) \cdot w.
$$
Corollary 6 implies that $g$ must divide $k+2$, $k-1$ or 1 and $g$
must divide $3k+1$ or $4k+2$. The first implies  that $k=1$ or
$k=g-2$. Since odd $g \ge 5$ cannot divide 4 or 6, we must have $k =
g-2$, so $g$ divides $3g-5$ or $4g-6$. This implies that $g=5$, so
$k=3$ and  $T = T_{2,17,25}$. However, $T$ contains the interior point
$(1,5,7)$, as well as  $w = (1,7,10)$, so is not a 1-point lattice
tetrahedron. 

Suppose $\ell = 3$, so $d_1=d_2=1$, $d_3=2$, $d_4=3$ and $N=7$.  It
follows from \eqref{E:v} that $3\ | \ 7e-2$, hence $e = 3k + 2$ and
$$
w = (1,3k+2,3g), \qquad v_4 = (2,7k+4,7g),
$$
with $0 \le k \le g-1$ and $\gcd(3k+2,3g) = \gcd(7k+4,7g) =
\gcd(7k+5,7g) = 1$. If $g=3$, the gcd conditions rule out
$k=1,2$, so $k=0$. But $T = T_{2,4,21}$ contains $(1,1,5)$ and
$(1,2,10)$ as well as $w = (1,2,9)$, and so is not a 1-point lattice
tetrahedron. Otherwise, $g \ge 5$, and we once again compute
$BC_{T_j}(0,0,1)$ for $j = 2,3$:
$$
(0,0,1) = \left( \frac{g+k+1}g\right) \cdot v_1 + \left(
\frac{-k}g\right) \cdot v_3 + \left( \frac 
1g\right) \cdot v_4 + \left( \frac{-2}g\right) \cdot w;
$$
$$
(0,0,1) =  \left( \frac{2g+3k+2}{2g}\right) \cdot v_1 + \left(
\frac{k}{2g}\right) \cdot v_2 +\left(  \frac 
{3k+2}{2g}\right) \cdot v_4 +\left(  \frac{-7k-4}{2g}\right) \cdot w.
$$
Again, $g$ divides 1, $k-1$ or $k+2$, so $k = 1$ or $k = g-2$, from
the first equation. If $k=1$, then $2g$ divides 6 or 10, so $g=5$ and
$k=3$, but $gcd(7\cdot3+4,7\cdot 5) > 1$. If $k=g-2$, then $2g$
divides 6 or 8, which is impossible.

Suppose $\ell = 4$, so $d_1=2$, $d_2=1$, $d_3=3$, $d_4=5$ and $N=11$.  It
follows from \eqref{E:v} that $5\ | \ 11e-3$, hence $e = 5k + 3$ and
$$
w = (1,5k+3,5g), \qquad v_4 = (2,11k+6,11g),
$$
with $0 \le k \le g-1$ and $\gcd(5k+3,5g) = \gcd(11k+6,11g) =
\gcd(11k+7,11g) = 1$. If $g=3$, the gcd conditions rule out
$k=0,1$, so $k=2$. But $T = T_{2,28,33}$ contains  $(1,6,7)$ 
and $(1,7,8)$ as well as $w = (1,13,15)$ and so is not a 1-point lattice
tetrahedron. Otherwise, $g \ge 5$, and we once again compute
$BC_{T_j}(0,0,1)$ for $j = 2,3$:
$$
(0,0,1) = \left( \frac{g+k+1}g\right) \cdot v_1 + \left(
\frac{-k}g\right) \cdot v_3 + \left( \frac 
1g\right) \cdot v_4 + \left( \frac{-2}g\right) \cdot w;
$$
$$
(0,0,1) =  \left( \frac{3g+5k+3}{3g}\right) \cdot v_1 + \left(
\frac{k}{3g}\right) \cdot v_2 +\left(  \frac 
{5k+3}{3g}\right) \cdot v_4 +\left(  \frac{-11k-6}g\right) \cdot w.
$$
Again, $g$ divides 1, $k-1$ or $k+2$, so $k = 1$ or $k = g-2$, from
the first equation. If $k=1$, then by the second equation, $3g$ divides 9
or 16, neither one of which is possible. If $k=g-2$, then $3g$
divides $6g-9$ or $10g-14$, which are also both impossible.

Suppose $\ell = 5$, so $d_1=1$, $d_2=3$, $d_3=4$, $d_4=5$ and $N=13$.  It
follows from \eqref{E:v} that $5\ | \ 13e-4$, hence once again $e =
5k+3$ and 
$$
w = (1,5k+3,5g), \qquad v_4 = (2,13k+7,13g),
$$
with $0 \le k \le g-1$ and $\gcd(5k+3,5g) = \gcd(13k+7,13g) =
\gcd(13k+8,13g) = 1$. If $g=3$, the gcd conditions rule out
$k=0,2,1$, in that order, so $g \ge 5$. As before,
$$
(0,0,1) = \left( \frac{3g+3k+2}{3g}\right) \cdot v_1 + \left(
\frac{-3k-1}{3g}\right) \cdot v_3 + \left( \frac 
1{3g}\right) \cdot v_4 + \left( \frac{-2}{3g}\right) \cdot w;
$$
$$
(0,0,1) =  \left( \frac{4g+5k+3}{4g}\right) \cdot v_1 + \left(
\frac{3k+1}{4g}\right) \cdot v_2 +\left(  \frac 
{5k+3}{4g}\right) \cdot v_4 +\left(  \frac{-13k-7}{4g}\right) \cdot w.
$$
From the first equation, $3g$ divides 1, $3k$ or $3k+3$, so $k=0$ or
$k = g-1$. If $k=0$, then the second equation implies that $4g$ divides
4 or 6. If $k = g-1$, then the second equation implies that $4g$
divides $8g-4$ or $10g-4$. None of these is possible. 

Suppose $\ell = 6$, so $d_1= 2$, $d_2=3$, $d_3=5$, $d_4=7$ and $N=17$.  It
follows from \eqref{E:v} that $7\ | \ 17e-5$, hence $e = 7k + 4$ and
$$
w = (1,7k+4,7g), \qquad v_4 = (2,17k+9,17g),
$$
with $0 \le k \le g-1$ and $\gcd(7k+4,7g) = \gcd(17k+9,17g) =
\gcd(17k+10,17g) = 1$. If $g=3$, the gcd conditions rule out
$k=2,0,1$, in that order, so  $g \ge 5$, and
$$
(0,0,1) = \left( \frac{3g+3k+2}{3g}\right) \cdot v_1 + \left(
\frac{-3k-1}{3g}\right) \cdot v_3 + \left( \frac 
1{3g}\right) \cdot v_4 + \left( \frac{-2}{3g}\right) \cdot w;
$$
$$
(0,0,1) =  \left( \frac{5g+7k+4}{5g}\right) \cdot v_1 + \left(
\frac{3k+1}{5g}\right) \cdot v_2 +\left(  \frac 
{7k+4}{5g}\right) \cdot v_4 +\left(  \frac{-17k-9}{5g}\right) \cdot w.
$$
From the first equation, $3g$ again divides 1, $3k$ or $3k+3$, so
$k=0$ or $k=g-1$. If $k=0$, the second equation implies that $5g$
divides 5 or 8. If $k = g-1$, then the second equation implies that
$5g$ divides $10g-5$ or
$14g-6$. Again, none of these is possible.

Finally, suppose $\ell = 7$, so $d_1=3$, $d_2=5$, $d_3=4$, $d_4=7$ and 
$N=19$.  It
follows from \eqref{E:v} that $7\ | \ 19e-4$, hence $e = 7k + 5$ and
$$
w = (1,7k+5,7g), \qquad v_4 = (2,19k+13,19g),
$$
with $0 \le k \le g-1$ and $\gcd(7k+5,7g) = \gcd(19k+13,19g) =
\gcd(19k+14,19g) = 1$. If $g=3$, the gcd conditions rule out
$k=1,2$, so $k=0$. But $T = T_{2,13,57}$, which contains at least
$(1,3,13)$ and
$(1,4,17)$ as well as $w = (1,5,21)$, and so is not a 1-point lattice
tetrahedron. Otherwise, $g \ge 5$, and, one last time:
$$
(0,0,1) = \left( \frac{5g+5k+4}{5g}\right) \cdot v_1 + \left(
\frac{-5k-3}{5g}\right) \cdot v_3 + \left( \frac 
1{5g}\right) \cdot v_4 + \left( \frac{-2}{5g}\right) \cdot w;
$$
$$
(0,0,1) =  \left( \frac{4g+7k+5}{4g}\right) \cdot v_1 + \left(
\frac{5k+3}{4g}\right) \cdot v_2 +\left(  \frac 
{7k+5}{4g}\right) \cdot v_4 +\left(  \frac{-19k-13}{4g}\right) \cdot w.
$$
The first equation implies that $5g$ divides one of 1, $5k+2$ or
$5k+5$, so $k = g-1$, and the second equation implies that $4g$
divides $12g-4$ or $14g-4$, neither of which is possible.

We have at long last completed a detailed proof of the following
theorem, which completes the classification of the 1-point tetrahedra.

\begin{theorem}
If $T$ is a 1-point tetrahedron with $BC_T(w) = \la^{(\ell)}$ and $g>1$,
then $\ell = 1$, $g=5$ and $T \approx T_{3,7,20} \approx T_{3,11,20}
\approx T_{7,11,20}$.
\end{theorem}

\section{Lattice widths and other questions}

If $S \in \rr^n$ is a lattice polytope and $u \in \zz^n$ then its
\emph{u-width} is defined to be $\max\{u\cdot x: x \in S\} -
\min\{u\cdot x: x \in S\}$, and its \emph{lattice width} is the
minimum of its $u$-widths, taken over $u \in \zz^n \setminus
0$. (Without loss of generality, we may always assume that the components
of $u$ have no common factor.)
Since $u\cdot x = (u M^{-1})\cdot(Mx)$, lattice width is preserved
by unimodular maps. If $S$ has lattice 
width $w$, then $S \cap \zz^n$ lies in $w+1$ consecutive ``lattice
hyperplanes'' $\pi_{j_0}, \dots \pi_{j_0+w}$, where  $\pi_j =
\{x: u\cdot x = j\}$. The Euclidean distance between $\pi_j$ and
$\pi_k$ is $|k-j|/|u|$, so a small geometric distance may correspond
to a large lattice width if $u$ has large components.

Theorem 5 shows that an empty lattice tetrahedron must have lattice
width 1. This does not hold for simplices in dimension $d \ge 4$; see
\cite{hz}. We present a possibly sporadic result for 1-point lattice
tetrahedra; it is proved using Theorem 7, rather than by an \emph{a
  priori} argument.

\begin{corollary} If $T$ is a 1-point lattice tetrahedron, then $T$
  has lattice width 2.
\end{corollary}
\begin{proof}
This is immediately true for any $T \approx T_{2,b,n}$, which is
contained in $0 \le x \le 2$; take $u = (1,0,0)$.
The two remaining cases are
$T_{3,3,4}$ and $T_{3,7,20}$, which are contained in $0 \le x+y-z \le
2$ and $0 \le 2x + 2y - z \le 2$, respectively. 
\end{proof}

More generally, suppose $T = T_{a,b,n}$ and $u = (r,s,t)$. Then we see
that the $u$-width of $T_{a,b,n}$ is equal to 
\begin{equation}\label{E:lw}
\max\{0,r,s,ar+bs+nt\} - \min\{0,r,s,ar+bs+nt\}.
\end{equation}
We present, without proof, the directions $\pm u$ in which the
1-point tetrahedra have width 2. Up to sign, the planes must be
$\{\pi_0, \pi_1, \pi_2\}$ if $v_1 = (0,0,0)$ is on an outer plane, or
$\{\pi_{-1}, \pi_0, \pi_1\}$ if $v_1$ is on the middle plane. This  gives us a
small finite set of $(r,s)$ to check. Let  $\ell_j = |T \cap \pi_j|$ and
let  $\ell_u(T):= (\ell_0, \ell_1, \ell_2)$  or
$\ell_u(T):= (\ell_{-1}, \ell_0,\ell_1)$, respectively. Since the
interior point must lie in the middle plane, the four possibilities
for  $\ell_u(T)$ are $(3,1,1)$ (or $(1,1,3), (2,1,2), (2,2,1)$ (or
(1,2,2)) and $(1,3,1)$. The first pair is impossible: suppose a face
of $T$ lies on a plane, then after a unimodular map we have $u = (0,0,1)$, 
the face can be placed on $z=0$, and has area $1/2$
by Pick's Theorem. By hypothesis, $T$ has altitude 2, and so volume
$\frac 26$, which is too small for a 1-point tetrahedron.
It turns out that $\ell_u(T)$ is {\it not} an invariant; several
of the smaller 1-point tetrahedra have different configurations in
different directions.

Somewhat surprisingly, $T_{3,3,4}$ has width two in nine directions:
$\ell_u(T)= (1,3,1)$ for $u = (1,0,-1), (0,1,-1), (1,-1,0), (2,1,-2),
(1,2,-2)$ and $(1,1,-1)$ and $\ell_u(T)= (2,1,2)$ for $u = (2,0,-1),
(0,2,-1)$ or $(2,2,-3)$. The next larger tetrahedron,  $T_{2,2,5}$, has
width two in six  directions: $\ell_u(T)= (1,3,1)$ for $u = (2,1,-1),
(1,2,-1)$ and $(1,-1,0)$ and  $\ell_u(T)= (2,2,1)$ for $u = (1,0,0),
(0,1,0)$ and $(1,1,-1)$. Next, $T_{2,4,7}$ has width two in four directions:
$\ell_u(T)= (1,3,1)$ for $u = (2,1,-1)$ and  $\ell_u(T)= (2,2,1)$ for
$u = (1,0,0), (0,2,-1)$ and $(1,1,-1)$. Each of $T_{2,6,11},
T_{2,7,13}$ and $T_{2,9,17}$ has width two in two directions:
$\ell_u(T)= (2,2,1)$ for $u = (1,0,0)$ and $u = (0,2,-1)$. Finally,
$T_{2,13,19}$ and $T_{3,7,20}$ each have width two in one direction: 
$\ell_u(T)= (2,2,1)$ for $u = (1,0,0)$ and $u = (2,2,-1)$
respectively. 

A check of 2-point lattice tetrahedra shows that most have lattice
width 2; however $T_{3,5,23}$ is one with lattice width 3.

We make some elementary remarks about the width
of lattice  $k$-tetrahedra for $k \ge 2$. As noted above,
$T_{3k,3k,3k+1}$ is a lattice $k$-tetrahedron; it also has width 2,
considering $u = (1,-1,0), (1,0,-1)$ or $(0,1,-1)$. Thus, width need
not go to infinity with $i(T)$. The following result applies to all
$T_{a,b,n}$, whether clean or not. 

\begin{theorem}
The lattice width of $T_{a,b,n}$ is $\le 2 \lceil n^{1/3} \rceil$.
\end{theorem}
\begin{proof}
This is a simple pigeonhole principle argument. Let $m = \lceil
n^{1/3} \rceil$ and consider $\{ ra + sb \mod n : 0 \le r,s \le
m\}$. There are $(m+1)^2$ residues, and so two must differ by at most
$(n-1)/(m+1)^2 < m$. Thus, $(r_1-r_2) a + (s_1 - s_2)b \equiv j \mod
n$ with $0 \le r_1, r_2 ,s_1, s_2, j \le m$. Now let $r = r_1-r_2$ and
$s=s_1-s_2$ and choose $t$ so that $r a + sb 
+ tn = j$ and let $u = (r,s,t)$.
 Since $-m \le r,s \le m$, \eqref{E:lw} implies the lattice width of
 $T$ in the $u$-direction is at most $2m$.
\end{proof}

The following example shows that the bound in Theorem 16 has the
correct order of magnitude. Let $T = T_{m,m^2,m^3+1}$. We first check
that this is clean. Clearly, $\gcd(m,m^3+1) =
\gcd(m^2,m^3+1) = 1$. Observe that $c(m,m^2,m^3+1) = 1 - m - m^2$ and 
let $h = \gcd(c,m^3+1)$. Since $c$ is odd, so is $h$, and since
$h$ is odd, and since $(m-1)c +(m^3+1) = 2m$, we have that $h \ |\
m$. But $h$ divides $m^3+1$, so $h=1$.
 
We claim that $T$ has lattice width $m$. This bound is achieved for $u
= (1,0,0)$ or $u = (1,m,-1)$. Suppose otherwise that
there exists $u = (r,s,t)$ so that
$$
\max\{0,r,s,mr+m^2s+(m^3+1)t\} - \min\{0,r,s,mr+m^2s+(m^3+1)t\} \le m-1.
$$
Then in particular, $|r|, |s| \le m-1$, and so $|mr + m^2s| \le
(m + m^2)(m-1) = m^3-m$. If $|t| \ge 1$, then
$$
|mr+m^2s+(m^3+1)t| \ge |(m^3+1)t| - (m^3-m) \ge m+1.
$$
If $t=0$, then $|mr + m^2s| = m|r + ms| \ge m$ unless $r+ms = 0$. In
this case, $m\ | \ r$ implies that $r=0$, so $s=0$ and $u = (0,0,0)$.
Therefore, the lattice width equals $m$. 

We make the following conjecture.
\begin{conjecture}
If $T$ is a $k$-point lattice tetrahedron, then its lattice width is 
$\le k+1$, and there is at least one interior lattice point on each of
the consecutive lattice planes in any minimal direction.
\end{conjecture}

If $T$ is a 2-point lattice  tetrahedron with interior points $w_1$
and $w_2$, then each interior point subdivides $T$ into 4
tetrahedra. The other interior point will be on an edge, a face or
interior to one of the subtetrahedra, and each case can occur.
For example, $T_{5,5,7}$ has two interior points: $w_1 = (1,1,1)$ and
$w_2 = (3,3,4)$, and $w_2 = \frac 12(v_4+w_1)$ is on an
edge, whereas  $w_1 = \frac 14(v_1+v_2+v_3+w_2)$ is interior.
Another 2-point lattice tetrahedron is $T_{5,5,8}$, with interior
points $w_1 = (1,1,1)$ and $w_2 = (2,2,3)$. Here, each is on a face
determined by the other: $w_1 = \frac13(v_2+v_3+w_2)$, $w_2 =\frac
13(v_1+v_4+w_1)$. It is possible for both to be interior; for example
if $T = T_{11,13,16}$ with $w_1 = (1,1,1)$ and $w_2 = (5,6,7)$, then
$$
BC_{T(v_1,v_2,v_3,w_2)}(w_1) = \left( \tfrac 37, \tfrac 27, \tfrac 17,
\tfrac 17 \right), \qquad
BC_{T(v_2,v_3,v_4,w_1)}(w_2) = \left( \tfrac 17, \tfrac 27, \tfrac 37,
\tfrac 17 \right).
$$

Another worthwhile project would be the classification of lattice
tetrahedra with one interior point and a positive number of boundary
points.
 Lemma 3 could be used in the special case that one of the four
faces has no non-vertex lattice points; the arguments of Theorem 4 can
be adapted to count the number of lattice points in $T_{a,b,n}$, when
\eqref{E:relpr} does not hold. It is not clear how to proceed if no face
is relatively empty. In view of White's Theorem, these questions
become considerably more difficult, even in four dimensions.

\bibliographystyle{amsalpha}

\end{document}